\documentclass[conference]{IEEEtran}
\IEEEoverridecommandlockouts
\usepackage{amsmath,amssymb,amsfonts,float}
\usepackage{algorithmic}
\usepackage{graphicx}
\usepackage{textcomp}
\usepackage{xcolor}
\usepackage{mathtools}
\usepackage{hyperref}
\usepackage{amsmath}
\usepackage{multirow}
\usepackage{booktabs}
\usepackage{threeparttable}
\usepackage{nccmath}
\begin{document}

\title{\huge Optimizing Demand-Responsive Paratransit Operations: \quad \quad A Mixed Integer Programming Approach
\thanks{This work was supported by the U.S. Department of Transportation through the Southeastern Transportation Research, Innovation, Development and Education (STRIDE) Region 4 University Transportation Center (Grant No. 69A3551747104).\\
\indent *Corresponding author.
}
}

\author{\IEEEauthorblockN{Xiaojian Zhang}
\IEEEauthorblockA{\textit{Dept. of Civil and Coastal Engineering} \\
\textit{University of Florida}\\
Gainesville, FL \\
xiaojianzhang@ufl.edu}\\
\IEEEauthorblockN{Noreen McDonald}
\IEEEauthorblockA{\textit{Dept. of City and Regional Planning} \\
\textit{UNC Chapel Hill}\\
Chapel Hill, NC \\
noreen@unc.edu}
\and
\IEEEauthorblockN{Yu Yang}
\IEEEauthorblockA{\textit{Dept. of Industrial and Systems Engineering} \\
\textit{University of Florida}\\
Gainesville, FL \\
yu.yang@ise.ufl.edu}\\
\IEEEauthorblockN{Xilei Zhao$^*$}
\IEEEauthorblockA{\textit{Dept. of Civil and Coastal Engineering} \\
\textit{University of Florida}\\
Gainesville, FL \\
xilei.zhao@essie.ufl.edu}
\and
\IEEEauthorblockN{Abigail L. Cochran}
\IEEEauthorblockA{\textit{Dept. of City and Regional Planning} \\
\textit{UNC Chapel Hill}\\
Chapel Hill, NC \\
acochran@unc.edu}

}

\maketitle

\begin{abstract}
The traditional demand-responsive paratransit system plays an important role in  connecting people to health care, particularly those who are carless, low-income, senior, underinsured/uninsured, or who have a disability. However, the existing paratransit system usually has low service quality due to long waiting times,  low  operation  frequencies,  and  high  costs.  In  order  to improve the service quality, we propose to design a new demand-responsive paratransit system that offers public, Uber-like options for non-emergency medical transportation. Mixed integer programming models are thus developed to optimize the system operations with the objectives of minimizing user waiting times for riders as well as operating costs for operators. The results produced in this paper will assist local departments of transportation and transit agencies as they consider operational strategies to meet non-emergency medical transport needs.
\end{abstract}


\section{Introduction}

There exist significant transportation barriers to health care facilities (e.g., hospitals, dialysis centers, and urgent care facilities) in the United States. A recent study found that, in 2017, 5.8 million Americans experienced delay in  non-emergency medical care due to a lack of transportation means \cite{wolfe2020transportation}. These people often have older age, lower incomes, disabilities, no access to personal vehicles, and/or limited or even no health insurance \cite{powers2016nonemergency,wolfe2020innovative}. Although the traditional demand-responsive paratransit systems may provide these people with access to health care facilities, the existing paratransit systems have suffered from long waiting times, low operation frequencies, and high costs \cite{kaufman2016intelligent}. 

In recent years, ridesourcing companies, such as Uber and Lyft, have emerged as important providers of non-emergency medical transportation services \cite{powers2016nonemergency,uberhealth}. In addition, health care providers are exploring the possibilities of using ridesourcing services to transport patients to and from medical appointments \cite{wolfe2020innovative}. Therefore, to meet the challenges of the changing market and ridership decline, many transit agencies are developing public, Uber-like options for non-emergency medical transportation, where customers can schedule a round trip from their home to a health care facility with preferred pick-up and drop-off times minutes prior to an appointment, or days in advance. However, how to design an efficient and economical paratransit system to provide such demand-responsive service remains largely unsolved.

One should note that the underlying problem, i.e., deploying the fleet to meet the non-emergency medical transportation needs in the most possibly efficient way, differs much from the problem faced by Uber and Lyft, where they try to maximize the profits by serving as many customers as possible in a timely manner. More specifically, Uber and Lyft have a large number of drivers that are scattered in the urban area. As a result, it is generally very likely to find available drivers that are close to customers. However, in our setting, we only a small fleet of vehicles to serve customers from a relatively large, rural area. In addition, non-emergency medical trips have much tighter time-window constraints for drop-offs since a late arrival to a health care facility is likely to result in a void trip and rescheduled appointment. More importantly, as paratransit is a service partially supported by the government, social equity plays an indispensable role. That means trips with a significantly longer travel distance have to be accepted by the operators at an affordable price, while commercial ridesourcing companies seldomly sacrifice their profits for equity purposes. To this end, substantial needs emerge for new models that are able to deploy the paratransit fleet efficiently and economically to meet the health care travel demands.

The essence of this problem is the trade-off between operating and waiting times. On the one hand, operators seek to minimize the total trip length, i.e., the total number of drivers' working hours, which affects the operating costs directly.  On the other hand, from the  customers'  perspective,  the  total waiting time, i.e., the difference between the actual and scheduled drop-off and pick-up times, is of great importance. Nevertheless, those two objectives cannot be achieved simultaneously in a straightforward way. If unlimited resources are available, we can simply dedicate a vehicle to transport each individual customer such that no waiting is ever needed, which inevitably drives up operating costs. However, with a fixed budget, we have to design a system that allows ridesharing to ensure each customer is picked up and dropped off as required while shortening the trip length as much as possible. The routing aspect of ridesharing has  to be taken into consideration, which complicates the problem drastically.

Almost all routing problems are notoriously difficult and mixed integer programming is one of the few exact solution frameworks that can yield high-quality solutions in a reasonable amount of time for small- to medium-sized problems.  
Actually, a very similar problem called the ``dial-a-ride problem'' (DARP) has been studied for a long time in the operations research community, whose  objective is to minimize the total traveling distance. Multiple models and algorithms have been proposed to tackle this problem, and \cite{cordeau2003dial} and \cite{cordeau2007dial} serve as good reviews on this topic. Traditional models for DARP try to minimize the total traveling distance, which is not the primary objective for paratransit services. Although a bunch of efficient heuristics have been developed to yield high-quality solutions (see \cite{cordeau2003tabu, attanasio2004parallel}) for DARP, we use an exact solution method in view of the moderate problem size and the needs for optimal solutions. More precisely, we approach this problem via the 3-index MIP model proposed in \cite{cordeau2006branch} with some modifications tailed to our problem of interest, which will be elaborated in Section \ref{Model}.

The main contributions of this study are summarized as follows:
\begin{itemize}
    \item We design a new, Uber-like paratransit service system from the operator's and the user's perspectives, respectively. The resulting problem is modeled as a mixed integer program (MIP).
    
    \item We propose two new objective functions in the two aforementioned situations, and add some application-specific constraints on top of the DARP model to accelerate the solution. 
    
    \item We use  a real-world data set as a case study and demonstrate that our approach is able to significantly improve the efficiency of the system.
\end{itemize}

The remainder of this paper is organized as follows: Section \ref{System Design and Basic Assumptions} is devoted to the problem description and  some basic assumptions. In Section \ref{Model Development}, we present the proposed two MIP models, i.e., the user model and the operator model. In Section \ref{Computational Experiments}, we explain the data set used for our case study and report the computational results. Section \ref{Conclusion} concludes our paper with pros and cons, and outlines future research directions.

\section{Problem Description and Assumptions}
\label{System Design and Basic Assumptions}
\subsection{Problem Description}
\label{System Description and Operating Policy}
\begin{figure}[htbp]
\centering
\includegraphics[width=1\linewidth]{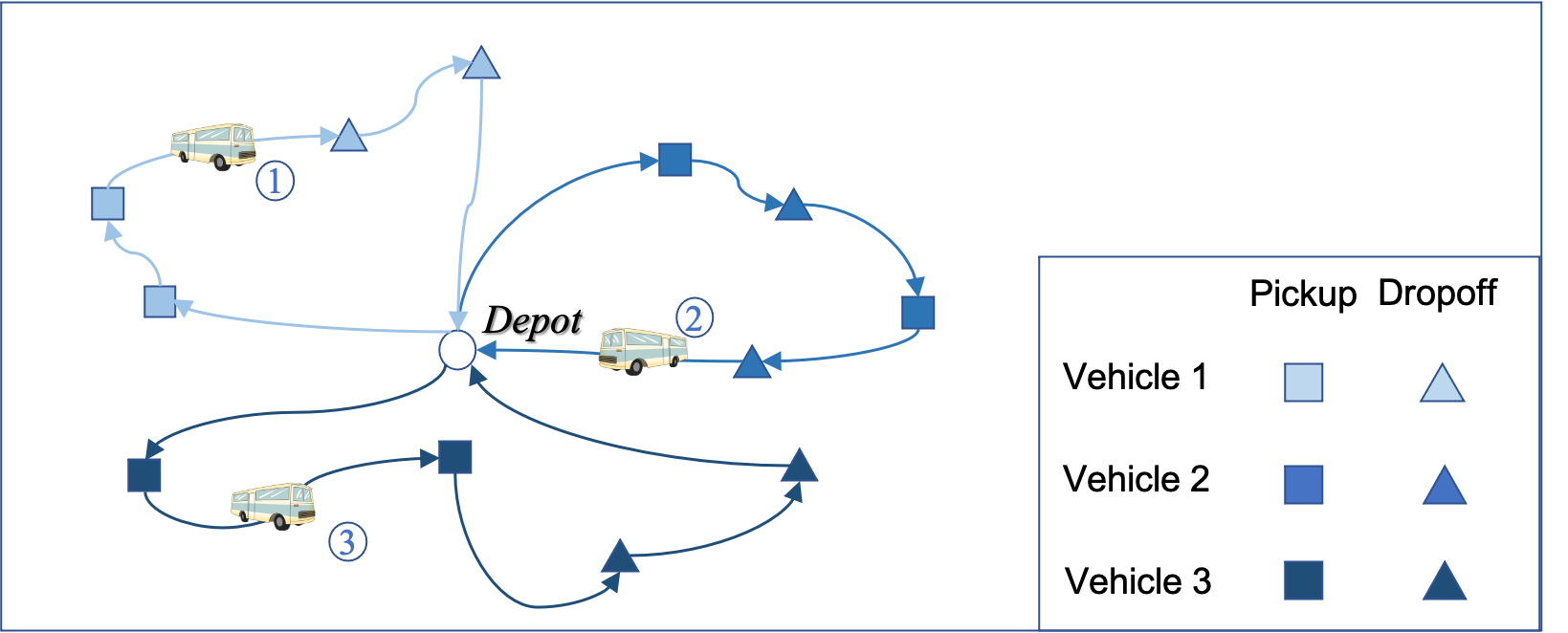}
\caption{Operating Policy}
\label{System Design}
\end{figure}

In this paper, we develop a Uber-like, door-to-door service system for providing non-emergency medical transportation. As shown in Fig. \ref{System Design}, each vehicle departs the depot to pick up and drop off a set of customers as required by their scheduled times and locations, and then returns to the depot after serving all the designated requests. The departure and return times for each vehicle are not necessarily the same. Usually, there will be an associated service time at each pick-up or drop-off location for boarding or alighting. All customers need to book their trips by calling or making a request online in advance, usually at least an hour before their desired pick-up time. Working hours are divided into intervals of equal length. Before the start of each working time interval, a group of vehicles are selected from the fleet to form a group to serve the requests within this time interval. The size of the group is dependent on the number of requests  received in this interval.


\subsection{Assumptions}
The assumptions made for modeling this problem are listed as follows.
\begin{itemize}
    \item The demand-responsive services are only available by booking at least one hour in advance, so the optimization can be done offline;
    \item A \(v_i\)-minute difference between the actual and scheduled pick-up/drop-off times at the location $i$ is allowed;
    \item All customers have no attendants. In other words, only one customer is served for each order;
    \item All vehicles meet the requirements of the Americans with Disabilities Act (ADA).
\end{itemize}

\section{Models}
\label{Model Development}
To provide a broader view for transportation and health policymakers, we analyze this problem from both the operator's and user's perspectives, which results in two models:  the \textbf{Operator Model (OM)} and the \textbf{User Model (UM)}. We first introduce the notation used throughout the rest of the paper.

\subsection{Notation}\label{Model}
Let $n$ denote the number of customers (orders) received within the time interval of interest. The model is constructed on a directed graph \(G = (N, A)\) with the node set \(N := \{0,1,...,2n+1 \}\) and the arc set $A$. Nodes \(0\) and \(2n+1\) represent the origin and destination depots, and subsets \(P = \{1,2,...,n \}\) and \(D = \{n,n+1,...,2n \}\) contain pick-up and drop-off nodes, respectively. Let $N^0 = N\setminus \{0, 2n+1\}$, then $A:=\{(i,j): \forall i, j\in N^0\}\cup \{(0,j):  \forall j\in P\}\cup \{(i,2n+1):  \forall i\in D\}$. Let \(K = \{1,2,...,p\}\) be the index set of the vehicles and vehicle $k$ has capacity \(C_k\). Each node  \(i \in N\) has a load \(q_i\), which is equal to 1 if $i\in P$, $-1$ if $i\in D$, and $0$ otherwise. Let \(d_i\) be the corresponding non-negative service time for $i\in N$. A time window \([e_i, l_i]\) is enforced for each node \(i \in N\) to make sure a vehicle will arrive within this time interval and the travel time is \(t_{i,j}\) for each arc $(i,j)\in A$.

\subsection{Decision Variables}
\begin{itemize}
\setlength\itemsep{-0.5em}
\item \(x_{ijk}\) (binary): equals 1 if vehicle \(k \in K\) uses arc $(i,j) \in A$, otherwise 0; \\
\item\(B_i\) (continuous): the time when vehicle \(k \in K\) arrives at node $i \in N$;          \\
\item\(Q_i\) (continuous): Number of customers on vehicle \(k \in K\) at node $i \in N$. Note the value should be integral, but it suffices to declare it to be continuous due to the model structure;;                     \\
\item\(y_i\)  (binary): indicator for potential waiting at node $i \in H$; where H is the set of all drop-off nodes of inbound trips and pick-up nodes of outbound trips;                  \\
\item\(z\) (continuous): objective to be optimized. \\
\end{itemize}

It should be noted that the first three variables, i.e., \(x_{ijk}\), \(B_i\), and \(Q_i\), are decision variables for OM, while all these five variables are decision variables for UM.

\subsection{Operator Model}
\label{OM}
From the operator's perspective, the goal is to serve the customers in a most cost-effective way while ensuring that each of them can  arrive  on time for the appointment and can be picked up from the health care facility no later than the scheduled time. Delayed or advanced pick-up's (drop-off's) from (at) home will be acceptable. A trip from home to a health care facility is called an \textbf{inbound} trip while one going back home is an \textbf{outbound} trip. In this model,  for a node $i$ that is the drop-off node of an inbound trip or a pick-up node of an outbound trip, the \(l_i\) is set to the scheduled drop-off/pick-up time. For the remaining nodes in $N^0$, the \(l_i\) is set to $L_i/2$ after the scheduled time, where $L_i$ is a predetermined number denoting the length of the time window. The earliest arrival time at a node $i$, denoted by $e_i$, is set accordingly as $L_i/2$ before the scheduled time to ensure that the length of the time window is equal to $L_i$.

In general, operators are more interested in the total operating costs, so the objective is set to minimizing the total operating time, $T$, of all vehicles as computed by \eqref{obj}.
Constraints \eqref{constr1} and \eqref{constr2} collectively ensure that every customer is visited only once and that the pick-up and drop-off nodes are visited by the same vehicle. Constraints \eqref{constr3} to \eqref{constr5} are  used to guarantee that each vehicle starts at the initial depot and finishes at the final depot. In scenarios where some of the vehicles are not used,  the vehicles leave the initial depot \(0\) and travel directly to  the final depot \(2n+1\) with 0 contribution to the objective value. 
Constraints \eqref{constr6} to \eqref{constr11} model the load and time relationships between successive nodes, where  \(M_1\) and \(M_2\) are  two sufficiently large constants that ensure the validity.
For example, if \(\sum_{k\in K} x_{i j k} = 1\), constraint \eqref{constr6} implies  a vehicle cannot arrive at node \(j\)  earlier than \(B_{i} + t_{ij} + d_i\) if it travels from node \(i\) to node \(j\). On the other hand, if \(x_{i j k} = 0\), \eqref{constr6} does not enforce any restriction. Constraint \eqref{constr12} ensures that each customer \(i\) will  be picked up before dropped off. Constraint \eqref{constr13} guarantees the each node \(i\)  is visited within a specific time window. Inequality \eqref{constr14} imposes the capacity constraint.

{
\footnotesize
\begin{align}
&\min \quad \sum_{k \in K} B_{2n+1, k} - B_{0k}\label{obj}\\
& \textit{subject to} \nonumber \\
&   \sum_{k \in K} \sum_{j \in N} x_{ijk}=1, \quad \forall i \in P, \label{constr1}\\
    &\sum_{j \in N} x_{ijk}-\sum_{j \in N} x_{n+i,jk}=0, \quad \forall i \in P, k \in K,\label{constr2}\\
    &\sum_{j \in N}  x_{0jk}=1, \quad \forall k \in K, \label{constr3}\\
    &\sum_{i \in N}  x_{i,2n+1,k}=1, \quad \forall k \in K,\label{constr4}\\
    &\sum_{i \in N} x_{jik}-\sum_{j \in N} x_{ijk}=0, \quad \forall i \in P \cup D, k \in K,\label{constr5}\\
    & B_{j} \geq B_{i}+t_{i j}+d_{i}-M_1\left(1-\sum_{k\in K}x_{i j k}\right), \forall i, j\in N^0, i\neq j,\label{constr6} \\
    & B_{2n+1,k} \geq B_{i}+t_{i,2n+1}+d_{i}-M_1\left(1- x_{i, 2n+1, k}\right),\forall i \in N^0, k\in K,\label{constr7} \\
 &B_{j} \geq B_{0 k}+t_{0 j}+d_{i}-M_1\left(1-x_{0 j k}\right), \forall j \in N^0 , k \in K, \label{constr8}\\
     & Q_{j} \geq Q_{i}+q_{j}-M_2\left(1-\sum_{k\in K}x_{i j k}\right), \quad  \forall i, j\in N^0, i\neq j, \label{constr9}\\
    & Q_{2n+1, k} \geq Q_{i} -M_2\left(1-x_{i,2n+1, k}\right), \quad \forall i \in N^0, k\in K, \label{constr10}\\
        & Q_{j} \geq q_j + M_2\left(1-x_{0jk}\right), \quad \forall i \in N, k\in K,\label{constr11}\\
    &  B_{i} + t_{i, n+i} + d_{i} \leq B_{n+i} \quad \forall i \in P,\label{constr12}\\
    &e_i \leq B_{i} \leq l_i \quad  \forall i \in N,\label{constr13}\\
    &\max \left\{0, q_{i}\right\} \leqslant Q_{i} \leqslant \min \left\{Q_{k}, Q_{k}+q_{i}\right\} \quad \forall i \in N, k \in K,\label{constr14}\\
    &x_{i j k} \in \left\{{0, 1}\right\} \quad \forall i \in N , j \in N , k \in K.\nonumber
\end{align}
}


In addition, to tighten the LP relaxation and accelerate the computation,  we also include the following two sets of constraints whose validity is straightforward.
\begin{align}
    & B_{2n+1,k} \geq B_{0k}, \quad \forall k \in K.
\end{align}

\subsection{User Model}

User experience is mostly dependent on the difference between the scheduled and actual drop-off times of inbound trips and pick-up times of outbound trips.
Thus, from the user's perspective, the goal is to minimize the sum of those differences. In reality, late arrival is less favorable than early arrival, especially, when significant delay occurs. To adjust the model away from excess lateness, a uniformly large penalty is incurred when the actual pick-up time is delayed more than a threshold $T$. To model such situations, we introduce a binary variable \(y_i\) as an indicator that takes a value of 1 if  the lateness is more than $T$, and equals to 0 otherwise. Let \(s_i\) be the scheduled time at node \(i\), and $\beta$ and $M_3$  be two large constants. Then the following constraint models the aforementioned situations.


\begin{equation}
    B_{i} - s_{i} \leq T(1-y_{i}) + M_3y_i, \quad\forall i\in H,\\
\end{equation}
The validity is due to the fact that $y_i$ is forced to be 1 when $B_i-s_i > T$, while it can either be 0 or 1 if $B_i-s_i \leq T$. 
 
Our objective is computed as the sum of the time difference or the potential penalty for excess delay, which can be modeled as follows.
\begin{align}
    &\min && z\nonumber \\
    &s.t. && z \geqslant \max \left\{\beta\sum_{i \in H} y_i, \sum_{i\in H} \lvert B_{i} - s_i \rvert \right\} \label{obj-user}
\end{align}
Note \eqref{obj-user} is  non-linear, which can be split into three linear constraints \eqref{user1} to \eqref{user3}:
\begin{align}
    &z \geqslant \beta\sum_{i \in H} y_i,\label{user1}\\
    &z \geqslant \sum_{i \in H} (B_{i} - s_i),\label{user2}\\
    &z \geqslant \sum_{i \in H}(s_i - B_{i}).\label{user3}
\end{align}
The complete user model will also need the constraints \eqref{constr1} to \eqref{constr14}. It should be noted that for $i\in H$,  \(e_i\) and \(l_i\) are set to \(0\) and \(1440\), respectively. For other nodes, the \(e_i\) and \(l_i\) are set in the  same way as in OM.

\subsection{Discussion}
The validity of the load and time constraints \eqref{constr6} to  \eqref{constr11} is ensured by  sufficiently large constants \(M_1\), \(M_2\) and \(M_3\). However, the larger these constants are, the looser the lower bound (the optimal values of the LP relaxation) tend to be. Thus, we would like to pick the smallest valid constants. In view of \(M_1 \geq \max\left \{B_{i}-B_{j}+t_{i,j}+d_{i}\right \}\), \(M_2 \geq \max\left\{Q_i - Q_j + q_i\right \}\) and \(M_3 \geq \max\left\{B_{i} - s_i\right\}\), we set \(M_1\)  to \(\max\left\{l_i\right\}-\min\left\{e_i\right\}+\max\left\{t_{ij}\right\}+\max\left\{d_i\right\}\), \(M_2\) to the maximum vehicle capacity, and  \(M_3\)  to  \(\max\left\{l_i\right\} - \min \left\{s_i\right\}\).

\begin{figure}[htbp]
\centering
\includegraphics[width=1\linewidth]{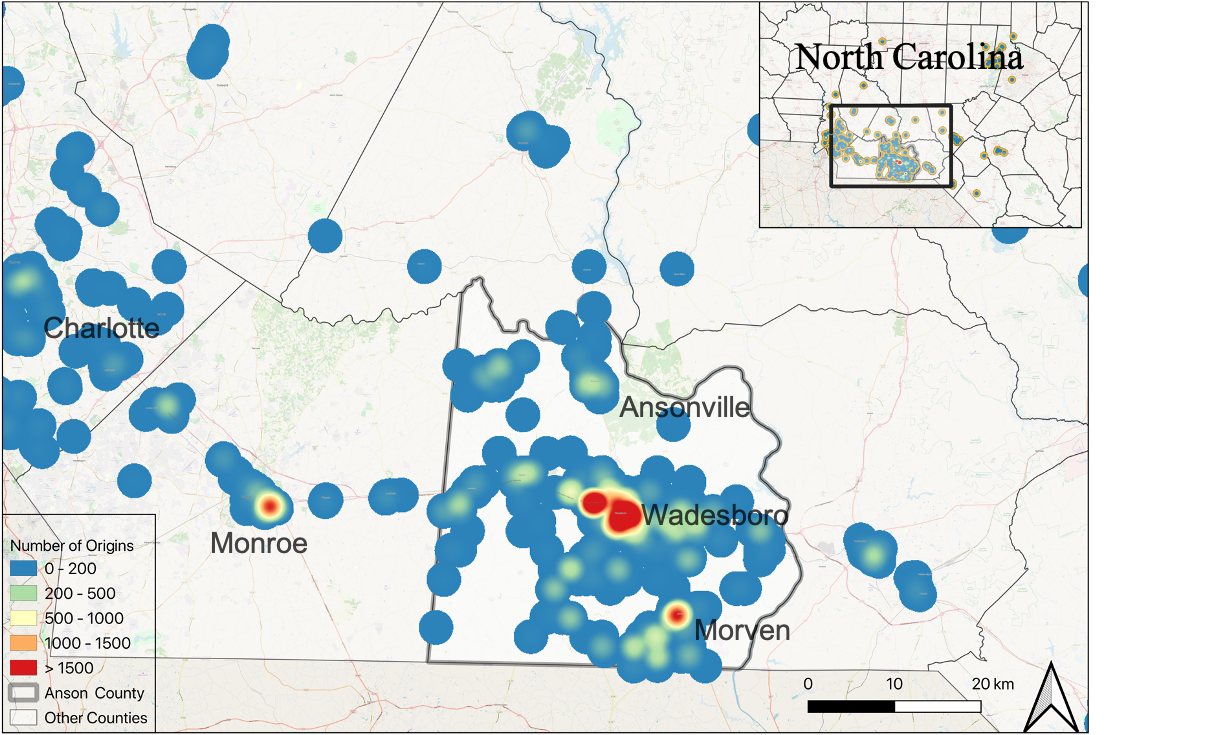}
\caption{Heatmap of trip-generation distribution}
\label{Heatmap of trip-generation distribution}
\end{figure}

\begin{table*}[h]
\centering
\caption{Model Summary}
\label{Number of Var Constrs and Run Time}
\resizebox{0.9\textwidth}{!}{%
\begin{tabular}{@{}lccccc|cccc@{}}
\toprule
\multirow{2}{*}{Period} & \multirow{2}{*}{\begin{tabular}[c]{@{}l@{}}\# of \\ orders\end{tabular}} & \multicolumn{4}{c|}{UM} & \multicolumn{4}{c}{OM} \\ \cmidrule(l){3-10} 
      &    & Vars & IntVars & Constrs & CPU (s) & Vars & IntVars & Constrs & CPU (s) \\ \midrule
5 am - 6 am   & 1  & 52   & 26      & 56     & 0.001    & 49   & 25      & 56      & 0.002    \\
6 am - 7 am   & 6  & 782  & 731     & 516    & 0.029    & 769  & 725     & 501    & 0.579    \\
7 am - 8 am   & 10 & 2086 & 2015    & 1172    & 0.638  & 2065 & 2005    & 1145    & 55.610   \\
8 am - 9 am   & 7  & 1048 & 992     & 656    & 0.052    & 1033 & 985     & 638    & 4.750    \\
9 am - 10 am  & 1  & 52   & 26      & 56     & 0.001    & 49   & 25      & 56      & 0.002    \\
10 am - 11 am & 8  & 1356 & 1294    & 815    & 0.995  & 1337 & 1285    & 791    & 6.176    \\
11 am - 12 pm & 5  & 556  & 510     & 392    & 0.062    & 545  & 505     & 380    & 0.188    \\
12 pm - 1 pm & 12 & 2978 & 2897    & 1596    & 75.030 & 2953 & 2885    & 1563    & 45.917   \\
1 pm - 2 pm & 3  & 224  & 188     & 192     & 0.006    & 217  & 185     & 186     & 0.016    \\
3 pm - 4 pm & 2  & 118  & 87      & 116     & 0.003    & 113  & 85      & 113     & 0.021    \\
4 pm - 5 pm & 3  & 224  & 188     & 192     & 0.009    & 217  & 185     & 186     & 0.092    \\ \bottomrule
\end{tabular}%
}
\end{table*}

\section{Computational Experiments}
\label{Computational Experiments}
In this section, we present the numerical results evaluating the performance of our proposed models, OM and UM, which are solved by the state-of-the-art MIP solver \textit{Gurobi} version 9.1.0. All experiments are implemented in Python and run on a workstation with Red Hat Enterprise Linux version 8.1, Intel(R) i9-9900K CPU @ 3.60GHz (8 physical cores, i.e.,  16 threads) and 64 GB of RAM. The time limit is set to an hour for each experiment, and all solver parameters are set to default.

\subsection{The Data}

We have access to the medical and nutritional purposes demand-responsive trips data collected by Anson County Transportation System (ACTS) in 2019. The data includes information about the scheduled and actual pick-up/drop-off timestamps and locations (i.e., latitude and longitude coordinates), appointment timestamps, odometer readings, cost billed (\$), dates, and use of mobility aids (e.g., wheelchairs). Since timestamps at pick-up and drop-off locations were manually recorded by drivers, errors were introduced inevitably. We treat trips with the same origin and destination, and those with travel distance less than 0 as outliers and remove them. We also remove incomplete data points, i.e., ones with missing values. After data cleaning, the total number of data points is 22,870 which consist of trips that took place on 261 different dates in 2019, and the average travel demand (origin-destination [OD] pairs) per day is 90. In addition, the trip starting times range from 3:00 am to 8:00 pm, and a major proportion (55.1\%) of the trips took place between 9:00 am and 1:00 pm. Spatially, as shown in Fig. \ref{Heatmap of trip-generation distribution}, most trips took place within Anson County, North Carolina, especially in Wadesboro and Morven, while a small fraction of the trips occurred outside Anson County, e.g., Monroe, Charlotte, and Durham. Moreover, around 65\% of the trips are short-to-medium-length trips with a travel distance less than 20 km (12.4 miles), while around 12\% are longer than 50 km (31.1 miles). Among all the trips, the shortest one is 0.13 km while the longest distance is 216 km.  We  use the \href{https://developers.google.com/maps/documentation/distance-matrix/overview}{Distance Matrix API} from \textit{Google Map} API to estimate the travel time and distance for each OD pair. 

\subsection{Experiment Description}
We used the trips that took place on January 3, 2019 as a case study. The total number of trips for this day is 58 (32 inbound trips to health care facilities and 26 outbound trips) and the scheduled times range from 5:00 am to 4:00 pm. We implement the OM and UM on an hourly basis,
which results in 11 different time intervals. The transit fleet size of ACTS is 14, and the capacity of each vehicle ranges from 7 to 18. As mentioned in Section \ref{System Description and Operating Policy}, we take the number of vehicles used within a time interval as an input parameter \(u\). In actual situations, the transit agency can flexibly select the number of vehicles served within the time interval \(I\). In this study, however, in order to consistently compare the results, we  set \(u\) to 5 uniformly and \(I\) to one hour. We also assume that all vehicles are identical with a maximum capacity of 7. In addition, we set the boarding time to be 7 minutes and alighting time to 5 minutes. According to ACTS, all vehicles were parked around their office location (i.e., the depot): 2485 US-74, Wadesboro, NC 28170. We also assume that all customers must allow a 30-minute time window for each  pick-up and drop-off. We set \(q_0 = q_{2n+1} = 0\), \(q_i + q_{n+i} = 0\), $\forall i \in P$ and \(d_0 = d_{2n+1} = 0\). The depot nodes, \(0\) and \(2n+1\), actually do need  a specific time window, but for consistency,  \((a_0, b_0)\) and \((a_{2n+1}, b_{2n+1})\) are both set to \((0, 1440)\). Note all numbers related to time windows have been converted into minutes.


\begin{table*}[htbp]
\centering
\caption{Results of the UM and OM.}\label{Results of the UM and OM}
\begin{threeparttable}
\resizebox{\textwidth}{!}{%
\begin{tabular}{@{}lccccccccc@{}}
\toprule
Period &
\# of orders &
  \begin{tabular}[c]{@{}l@{}} UM\_Raw (min)\end{tabular} &
  \begin{tabular}[c]{@{}l@{}}UM\_UM (min)\end{tabular} &
  \begin{tabular}[c]{@{}l@{}}UM\_OM (min)\end{tabular} &
  OM\_UM (min)&
  OM\_OM (min)&
  V\_UM &
  V\_OM \\ \midrule
5 am - 6 am   & 1  & 15   & 0 (15)     & 9 (6)       & 30   & 30   & 1  & 1  \\
6 am - 7 am  & 6  & 96  & 0 (96)   & 116 (--20)   & 547  & 408  & 3  & 2  \\
7 am - 8 am  & 10 & 190  & 0 (190)   & 187 (3)   & 1173 & 516  & 5  & 3  \\
8 am - 9 am  & 7  & 163  & 0 (163)    & 48 (115)     & 553  & 359  & 5  & 1  \\
9 am - 10 am & 1  & 10   & 0 (10)     & 30 (--20)     & 33   & 33   & 1  & 1  \\
10 am - 11 am & 8  & 68  & 8 (60)   & 145 (--77)    & 760  & 502  & 5  & 3  \\
11 am - 12 pm & 5  & 135  & 0 (135)    & 23 (112)   & 571  & 267  & 5  & 2  \\
12 pm - 1 pm & 12 & 926 & 45 (881)  & 126 (800)  & 956  & 736  & 5  & 5  \\
1 pm - 2 pm & 3  & 228  & 0 (228)    & 20 (208)    & 261  & 220  & 3  & 2  \\
3 pm - 4 pm & 2  & 83  & 0 (83)    & 38 (45)    & 65   & 61   & 2  & 1  \\
4 pm - 5 pm & 3  & 156  & 0 (156)    & 69 (87)   & 141  & 100  & 3  & 1  \\
\midrule
\qquad Total & 58 & 2070 & 53 (2017) & 811 (1259) & 5090 & 4124 & 38 & 23 \\ \bottomrule
\end{tabular}%
}
\end{threeparttable}
\end{table*}

\subsection{Results}

Table \ref{Number of Var Constrs and Run Time} presents
statistics about the models: Number of orders, number of variables, number of constraints, and solution time. 
All the cases can be solved to optimality within an hour. Table \ref{Results of the UM and OM} summarizes the results of the UM and OM. For a better understanding of the trade-off involved, we compute the UM objective (without the large penalty for excess lateness) using the solution yielded by the OM and vice versa. For convenience, we use A\_B to denote the value computed by the objective function of model A at the solution yielded by model B. Thus, UM\_Raw, UM\_UM, UM\_OM in  Table \ref{Results of the UM and OM} are the UM objective values ($\sum_{i\in H} \lvert B_{i} - s_{i}\rvert$) evaluated at the existing operational data, the solutions yielded by UM and OM, respectively. The number in each bracket represents the reduction compared to UM\_Raw. OM\_UM, OM\_OM are the OM objective values evaluated at the solutions yielded by UM and OM, respectively. Lastly, V\_UM and V\_OM are the number of vehicles actually used by UM and OM, respectively.

As shown in Table \ref{Number of Var Constrs and Run Time}, in many cases, UM has slightly more variables and constraints than OM. We observed that when the number of trips is smaller than 8, both models can be solved in seconds. However, as the number of trips increases, solving OM becomes more time-consuming, but still less than
1 minute. It is worth mentioning that while solving UM for most instances is efficient, for period 12 pm - 1pm, the CPU reaches more than 75 seconds. A possible explanation is that this period has the greatest number of orders; and compared to other periods, the spatial and temporal distributions of orders in this period are more uneven, which largely increases the CPU computing time.

The third column of Table \ref{Results of the UM and OM} shows the time difference of nodes in $H$ calculated by the raw data, which varies significantly across different time periods. For example, from 7 am to 8 am, the difference is 190 minutes for 10 trips, while from 12 pm to 1 pm, it is 926 minutes for 12 trips. In addition, based on the raw data, the average time difference for each trip is 35.7 min. In contrast, our proposed UM yields substantially better results where the average is reduced to 0.9 min for each trip. The UM can also improve this metric by around 97.4\%. As mentioned, we evaluate the time difference of the solution yielded by OM, which is shown in the  fifth column of Table \ref{Results of the UM and OM}. 
In addition, we observe that there are some instances (e.g., time period 6 am - 7 am) whose time difference of OM is worse than that of the benchmark, which is probably due to the modeling logic of OM described in Section \ref{OM}. More specifically, from the operator's perspective, the operating policy is to ensure all customers reach and leave hospitals on time while minimizing the total operating time. Hence, in order to lower the total operating time, fewer vehicles will be used and more ridesharing will occur, resulting in an  increase in customers' in-vehicle time and thus an increase in  the time difference.

According to the results shown in Table \ref{Results of the UM and OM}, it is clear that the total operating time of OM (4124 minutes) is less than that of UM (5090 minutes). However, this improvement is not very significant, which indicates that UM does not sacrifice much on the total vehicle operating time even though its objective is to minimize the total time difference. Since improving users experience is of great importance to transit agencies, it is advisable for key stakeholders like cities and transit agencies to leverage the proposed UM when planning new, demand-responsive paratransit systems. Moreover, another finding is that UM generally uses more vehicles than OM, which is reasonable since using fewer vehicles will reduce the total operating time and lead to more cost-effective operations of the paratransit system.


\section{Conclusion}
\label{Conclusion}

This paper designs a novel, Uber-like paratransit system from both the operator's and the user's perspectives, which is solved by using a MIP approach. Compared to the current paratransit service, the developed UM can not only considerably reduce the time difference between the actual and scheduled times (i.e., a 97.4\% reduction), but can also achieve a satisfying total operating time. In other words, the UM prioritizes the user experience without significantly sacrificing the operations efficiency. Cities, local departments of transportation, and transit agencies should consider adopting the UM when they develop new strategies to meet non-emergency medical transport needs. 

The paper has some limitations. For the computational experiments, we  only use the data from a single day. More data from multiple days should be used to test the proposed models. Furthermore, all vehicles are considered to be identical in the experiments. But in reality, the capacity of the vehicles may be different and only a fraction of the fleet are ADA accessible. Therefore, future work will be focused on producing more realistic results by taking these elements into consideration.

\section*{Acknowledgment}
The authors would like to thank Kai Monast for securing, cleaning, and explaining the data.

\bibliographystyle{IEEEtran}
\bibliography{reference}

\end{document}